\documentclass[12pt]{amsart}
\usepackage{amscd,amssymb,graphicx}
\newif\ifpdf\ifx\pdfoutput\undefined\pdffalse\else\pdftrue\fi
\ifpdf\usepackage[colorlinks,pagebackref]{hyperref} 
\else\usepackage[hypertex,pagebackref]{hyperref}\fi 
\def\psphere{N}
\renewcommand\deg{\mathop{\mathsf{deg}}}
\newcommand\Ker{\mathop{\mathsf{Ker}}}

\newcommand\GL{\mathop{\mathsf{GL}}}
\newcommand\PGL{\mathop{\mathsf{PGL}}}
\newcommand\PSL{\mathop{\mathsf{PSL}}}
\newcommand\Aff{\mathop{\mathsf{Aff}}}
\newcommand\tr{\mathop{\mathsf{tr}}}
\newcommand\SL{\mathop{\mathsf{SL}}}

\newcommand\Out{\mathop{\mathsf{Out}}}
\newcommand\Aut{\mathop{\mathsf{Aut}}}

\newcommand\C{{\mathbb C}}
\newcommand\R{{\mathbb R}}
\newcommand\Z{{\mathbb Z}}
\newcommand\A{{\mathbb A}}

\newcommand\pgltz{\PGL(2,\Z)}
\newcommand\homeo{\mathop{\mathsf{Homeo}}}
\newcommand\slt{\SL(2,\C)}
\newcommand\mm{\mathcal M}
\renewcommand\P{\mathbb{P}}
\newtheorem{theorem}{Theorem}[section]
\newtheorem{itheorem}{Theorem}
\newtheorem{iproposition}[itheorem]{Proposition}
\newtheorem*{theorem*}{Theorem}
\newtheorem{lemma}[theorem]{Lemma}

\newtheorem*{proposition*}{Proposition}

\newtheorem*{conjecture*}{Conjecture}

\theoremstyle{definition}

\newtheorem*{remark*}{Remark}
\newtheorem*{acknowledgments}{Acknowledgments}

\begin{document}

\title[Modular group action]
{Homological action of the modular group on some cubic moduli spaces}
\author{William M. Goldman}
\thanks{Research supported under NSF grant
DMS-0103889.}
\address{Department of Mathematics\\
The University of Maryland\\College Park, MD 20742}
\email{wmg@math.umd.edu}
\author{Walter D. Neumann}
\thanks{Research supported under NSF grants DMS-0083097 and DMS-0206464}
\address{Department of Mathematics\\Barnard College, 
Columbia University\\New York, NY 10027}
\email{neumann@math.columbia.edu}
\begin{abstract}
  We describe the action of the automorphism group of the complex
  cubic $x^2+y^2+z^2-xyz-2$ on the homology of its fibers. This action
  includes the action of the mapping class group of a punctured torus
  on the subvarieties of its $\SL(2,\C)$ character variety given by
  fixing the trace of the peripheral element (so-called ``relative
  character varieties''). This mapping class group is isomorphic to
  $\PGL(2,\Z)$. 

  We also describe the corresponding mapping class group
  action for the four-holed sphere and its relative $\SL(2,\C)$
  character varieties, which are fibers of deformations
  $x^2+y^2+z^2-xyz-2-Px-Qy-Rz$ of the above cubic. The $2$-congruence
  subgroup $\PGL(2,\Z)_{(2)}$ still acts on these cubics and is the
  full automorphism group when $P,Q,R$ are distinct.
\end{abstract} 
\maketitle
\tableofcontents

\section*{Introduction}

Several important moduli spaces are described by cubic surfaces in
affine space. The symmetries of the original moduli problem determines
a group of transformations of the moduli space which generates a
dynamical system. For two bounded surfaces $M$ ---
the one-holed torus and the four-holed sphere --- 
the relative $\slt$-character varieties
are described by cubic surfaces $\mm$ in $\C^3$ of the form
\begin{equation*}
x^2 + y^2 + z^2 - x y z = P x + Q y + R z + S 
\end{equation*}
upon which the mapping class group $\pi_0\homeo(M,\partial M)$
acts by polynomial automorphisms. For these surfaces
$\pi_0\homeo(M,\partial M)$ is essentially the modular
group $\pgltz$. We compute the homology of $\mm$ and its intersection
form, and the resulting action of $\pgltz$ on $H_*(\mm)$.

A {\em relative $\slt$-character variety\/} of a compact bounded surface 
$M$ classifies equivalence classes of representations 
\begin{equation*}
\rho:\pi_1(M)\longrightarrow\slt 
\end{equation*} 
with certain boundary conditions. Namely let $\partial_iM$ for
$i=1,\dots,n$ denote the boundary components of $M$ and choose
$\delta_i\in\pi_1(M)$ corresponding to generators of
$\pi_1(\partial_iM)$ respectively and and
$t=(t_1,\dots,t_n)\in\C^n$. Then the relative $\slt$-character variety
with boundary data $t$ is the moduli space of representations $\rho$
where
\begin{equation*}
\tr\big(\rho(\delta_i)\big) = t_i
\end{equation*}
for $i=1,\dots,N$. 

Suppose first that $M$ is a one-holed torus. Its fundamental group
$\pi$ is free of rank two, and admits a {\em redundant geometric\/} 
presentation 
\begin{equation*}
\langle X,Y, \delta_1 \mid \delta_1 = XYX^{-1}Y^{-1} \rangle
\end{equation*}
where $\delta_1$ corresponds to a generator of $\pi_1(\partial M)$.
Fricke \cite{fricke} proved that the character
variety is an affine space $\C^3$ with coordinates
\begin{align*}
x &= \tr(\rho(X)) \\ 
y &= \tr(\rho(Y)) \\ 
z &= \tr(\rho(Z)),\quad Z=XY\,.
\end{align*}
The resulting $\pi_0\homeo(M,\partial M)$-action on this $\C^3$ 
factors through an action of
$\PGL(2,\Z)$. Namely, by Dehn (unpublished) and 
Nielsen \cite{nielsen}, $\pi_0\homeo(M,\partial M)\cong\Out(\pi)\cong
\GL(2,\Z)$ via the action on the abelianization $\Z^2$
of $\pi=\pi_1(M)$. The center $\{\pm I\}$ of $\GL(2,\Z)$ is represented by the
automorphism $X\mapsto X^{-1}, Y\mapsto Y^{-1}, Z\mapsto
YZ^{-1}Y^{-1}$, which preserves traces and hence acts trivially on
$\C^3$, so the action can be considered as an action of
$\GL(2,\Z)/\{\pm I\}=\PGL(2,\Z)$.

The relative character varieties result by fixing the trace of 
$\rho(\delta_1)$, which is given by the polynomial
\begin{equation*}
\tr(\rho(XYX^{-1}Y^{-1})) = \kappa(\tr(\rho(X)),\tr(\rho(Y)),\tr(\rho(XY))
\end{equation*}
where 
\begin{equation*}
\kappa(x,y,z) := 
 x^2 + y^2 + z^2 - x y z - 2. 
\end{equation*}
Thus, the action of $\PGL(2,\Z)$ on $\C^3$ preserves this polynomial
$\kappa(x,y,z)$.

The full automorphism group $\Gamma=\Aut(\kappa)$ of this polynomial
$\kappa$ is generated by $\PGL(2,\Z)$ and the ``group of sign
changes,'' namely the Klein 4-group $\Sigma\cong C_2\times C_2$ that
acts on $\C^3$ by $(x,y,z)\mapsto (\pm x,\pm y,\pm z)$ with an even
number of minus signs.  In fact (see Horowitz \cite{horowitz})
$\Gamma$ is the split extension:
$$\Gamma=\PGL(2,\Z)\ltimes\Sigma.$$
The dynamic behavior of the
action of $\Gamma$ on the real character variety $\R^3\subset\C^3$ was
described in \cite{goldman}, but little is known about the full action
on $\C^3$. In this note we investigate the algebraic topology of this
action of $\Gamma$ (and hence of $\PGL(2,\Z)$) on the fibers of the
polynomial $\kappa$.

The results can be summarized:
\begin{itheorem}\label{th:i1}
The fibers $V_t:=\kappa^{-1}(t)$ of the polynomial 
$\kappa\colon \C^3\to \C$ have reduced homology only in dimension 2. 
There are two
special fibers: $V_2$ and $V_{-2}$ with homology $H_2(V_{2})\cong \Z$ and 
$H_2(V_{-2})\cong \Z^4$, respectively. 
The generic fiber $V_{t}$ with $t\ne\pm 2$ has homology 
\begin{equation*}
H_2(V_{t})= H_2(V_{-2})\oplus   H_2(V_{2})\cong \Z^5.  
\end{equation*}
The automorphism group $\Gamma$ acts on these
homology groups via a homomorphism to $S_4\times C_2$, with the
symmetric group $S_4$ acting by $(\rho_4)\oplus1$ on $\Z^4\oplus\Z$,
where $\rho_4$ is the standard permutation representation of $S_4$
on $\Z^4$, and the cyclic group $C_2$ acting by multiplication by
$\pm1$ on $\Z^5$.

The mapping class group $\Out(\pi)$, acting on $\C^3$ as the subgroup\break
$\PGL(2,\Z)$ of $\Gamma$, acts on the above homology groups as the
subgroup $S_3\times C_2$ of $S_4\times C_2$, where the $S_3$ is the
subgroup of $S_4$ that fixes one of the generators of $\Z^4$.
\end{itheorem}

The homomorphism of $\Gamma$ to $S_4\times C_2$ is  as
follows: The homomorphism to $C_2$ is the determinant map
$$\Gamma=\PGL(2,\Z)\ltimes\Sigma\to\PGL(2,\Z)
\stackrel{\text{\tiny det}}\longrightarrow C_2$$
The homomorphism to $S_4$ arises via reduction modulo $2$:
\begin{equation*}
\Gamma=\PGL(2,\Z)\ltimes\Sigma~
\longrightarrow~\PGL(2,\Z/2)\ltimes\Sigma=S_3\ltimes\Sigma 
\cong S_4.
\end{equation*}
Here $\PGL(2,\Z/2)\ltimes\Sigma \cong S_4$  because
$\PGL(2,\Z/2)\ltimes\Sigma$ is the group of affine automorphisms of an
affine plane $\A = \A^2(\Z/2)$ and every permutation of the
$4$-element set $\A$ is affine.  The subgroup $\Sigma$ 
is the group of translations, which acts simply transitively on $\A$.
The group $S_3$ that fixes a chosen origin in $\A$ and permutes the
other three points is the group of linear automorphisms
$\PGL(2,\Z/2)=\GL(2,\Z/2)$ (since $\Z/2$ contains only one nonzero
scalar,
\begin{equation*}
\PGL(2,\Z/2)=\GL(2,\Z/2)=\PSL(2,\Z/2)=\SL(2,Z/2)~). 
\end{equation*}

In terms of the punctured torus $M$,  $\A$ can be identified with the set 
of spin structures on a trivial vector bundle; 
specifically its vector space of translations identifies
with the cohomology group $H^1(M,\Z/2)$.

Note that the kernel of $\Gamma= \PGL(2,\Z)\ltimes\Sigma
\longrightarrow \PGL(\Z/2,2)\ltimes\Sigma=S_4$ is the congruence subgroup
$$\PGL(2,\Z)_{(2)}:=\left\{
\begin{pmatrix}
  a&b\\c&d
\end{pmatrix}
\in\PGL(2,\Z)~|~b\equiv c\equiv 0\text{ mod }2\right\}$$
of
$\PGL(2,\Z)$.
We have a commutative diagram with exact rows and columns
$$
\begin{CD}
&&&&  1&&1\\
&&&&@VVV @VVV\\
&&&& \Sigma @>=>> \Sigma\\
&&&&@VVV @VVV\\
1 @>>> \PGL(2,\Z)_{(2)} @>>> \Gamma @>>> S_4 @>>> 1\\
&& @VV=V @VVV @VVV\\
1 @>>> \PGL(2,\Z)_{(2)} @>>> \PGL(2,\Z)@>>> S_3 @>>> 1\\
&&&&@VVV @VVV\\
&&&&  1&&1\\
\end{CD}
$$
The right-hand column identifies with the short exact sequence
$$1\longrightarrow\Z/2\oplus\Z/2 \longrightarrow \Aff(2,\Z/2)
\longrightarrow \PGL(2,\Z/2) \longrightarrow1$$
All four short exact sequences of the diagram are split exact. In
particular,
\begin{equation}
\Gamma\cong\PGL(2,\Z)_{(2)}\rtimes S_4\,.
  \label{eq:1}
\end{equation}

We now describe the case of the four holed sphere $\psphere$. The
group $\Gamma$ arises in a new context in this discussion, as the
mapping class group $\pi_0\homeo(\psphere,\partial\psphere)$. The
relative $\SL(2,\C)$ character varieties for $\psphere$ are again certain
cubic surfaces. We first discuss the automorphisms of these cubic
surfaces.

\begin{itheorem}\label{th:i2}
  For any $P,Q,R\in\C$, the action of $\PGL(2,\Z)_{(2)}$ on $\C^3$
  that preserves the polynomial $\kappa$ deforms to an action that
  preserves the cubic polynomial
  $$\kappa_{P,Q,R}(x,y,z):=x^2 + y^2 + z^2 - x y z - P x - Q y - R
  z-2\,.$$
  Each nonsingular fiber $V$ of $\kappa_{P,Q,R}$ is
  diffeomorphic to a non-singular fiber of $\kappa$, so it still has
  homology $H_2(V) \cong \Z^5$, while a singular fiber has homology a
  quotient of this. The action of $\PGL(2,\Z)_{(2)}$ on the homology
  of fibers of $\kappa_{P,Q,R}$ is again by multiplication by $\pm1$
  via the determinant homomorphism $\PGL(2,\Z)_{(2)}\to C_2=\{\pm1\}$.

  $\PGL(2,\Z)_{(2)}$ is the full automorphism group of
  $\kappa_{P,Q,R}$ if $P,Q,R$ are distinct, while otherwise the full
  automorphism group is $\PGL(2,\Z)_{(2)}$ extended by an appropriate
  subgroup of $S_4$ (those permutations and sign changes of $x,y,z$ 
  that preserve $\kappa_{P,Q,R}$; --- if at
  most one of $P,Q,R$ is zero then there are only permutations).
\end{itheorem}
Each relative character variety for $\psphere$ is a fiber of some
cubic polynomial $\kappa_{P,Q,R}(x,y,z)$. The mapping class group
that fixes all four boundary components satisfies
$$\pi_0\homeo(\psphere,\partial_1 \psphere,\dots,\partial_4 \psphere) \cong
\PGL(2,\Z)_{(2)}\,,$$
and its action on the relative character variety is
by the action of the above theorem.

If the trace constraints at some boundary components of $M$ are
identical, then 
the relevant mapping class group is $\PGL(2,\Z)_{(2)}\rtimes K$ where
$K\subset S_4$ is the group of permutations that preserve trace
constraints. Its action on the relative character variety is by the
quotient $\PGL(2,\Z)_{(2)}\rtimes \pi(K)$ where $\pi\colon S_4\to
S_3$.  In particular, if the trace constraints at all four boundary
components of $\psphere$ are equal, then the relevant mapping class
group is a semidirect product $\PGL(2,\Z)_{(2)}\rtimes S_4$. In fact
it is the same semidirect product as that of equation (\ref{eq:1}),
so 
$$\pi_0\homeo(\psphere,\partial\psphere)\cong \Gamma\,.$$
But its
action on $\C^3$ in Theorem \ref{th:i2} is via its homomorphism to
$\PGL(2,\Z)$, with kernel $\Sigma$, while the action of $\Gamma$ on
$\C^3$ of Theorem \ref{th:i1} was faithful.

\vspace{3pt}We also determine the intersection form on homology of fibers:

\begin{iproposition}\label{prop:i1}
The intersection form on the homology
$H_2(V_t)$ of a general fiber of any one of the cubics
$\kappa_{P,Q,R}$ with respect to the basis used in the
above theorems is
$$\begin{pmatrix}
  -2&0&0&0&1\\
0&-2&0&0&1\\
0&0&-2&0&1\\
0&0&0&-2&1\\
1&1&1&1&-2
\end{pmatrix}\,.$$
\end{iproposition}

\begin{acknowledgments}
Goldman is grateful to J.\ Damon, L.\ Ein, H.\ King, A.\
Libgober and J.\ W.\ Wood for valuable conversations.
\end{acknowledgments}

\section{Vanishing cycles and homology}\label{sec:homology}

It is known (e.g., \cite{neumann-norbury}) that the homology of a
general fiber of a polynomial map $\kappa$ comes from vanishing cycles
of the singularities and vanishing cycles ``at infinity.'' 

Broughton \cite{broughton} shows that a polynomial has no vanishing
cycles at infinity if it is ``\emph{tame}'', that is,
$|d\kappa|$ is bounded away from zero outside some compact set. Our
polynomial $\kappa$ is tame, in fact:
\begin{lemma}\label{le:grad}
$|d\kappa|$ tends to infinity as
$|(x,y,z)|\to\infty$.   
\end{lemma}
\begin{proof}
\begin{equation*}
d\kappa = (2x-yz) dx + (2y-xz) dy + ( 2z-xy) dz =: a dx + b dy + c dz
\end{equation*}
We shall show that $\Vert(a,b,c)\Vert$ bounded implies 
$\Vert(x,y,z)\Vert$ bounded.  
Consider a sequence of points $(x,y,z)$ going to infinity for
which $(a,b,c)$ stays bounded.  If two coordinates, say $x$ and $y$,
stay bounded then $c=2z-xy$ tends to infinity, a contradiction.
So at least two coordinates, say $y$ and $z$, must approach
infinity in the sequence. But for $z\ne\pm2$ the equations $2x=a+yz$
and $2y=b+xz$ imply 
\begin{equation*}
y=\frac{2b+az}{4-z^2}, 
\end{equation*}
which approaches $0$ as 
$z\to\infty$, another contradiction.
\end{proof}
For a tame polynomial, the homology of the general fiber is given by
the vanishing cycles of singularities, so we need to find the
singularities of $\kappa$:  
\begin{equation*}
2x-yz = 2y-xz =  2z-xy =0 
\end{equation*}
implies the critical points of $\kappa$ are:
\begin{equation}\label{eq:cycles}
(2,-2,-2),~(-2,2,-2),~(-2,-2,2),~(2,2,2),~(0,0,0).
\end{equation}
The first four 
singular points are on the fiber $V_2$ while the singular point
$(0,0,0)$ is on $V_{-2}$.  For any one of the singular points
$(2\epsilon_1 ,2\epsilon_2 ,2\epsilon_1\epsilon_2)$ with $\epsilon_i=\pm1$ we
introduce local coordinates
\def\x{\hat x}
\def\y{\hat y}
\def\z{\hat z}
$$\x=\epsilon_1 x-2,\quad \y=\epsilon_2 y-2,\quad
\z=\epsilon_1\epsilon_2 z-2\, .$$
In these local coordinates we have
$$\kappa(x,y,z)=\x^2+\y^2 +\z^2 -2\x\y-2\x\z-2\y\z-\x\y\z.$$
The Hessian at
$(0,0,0)$ is $2I$ and at each of the other four singular points
(using the above local coordinates) it is
$$2\begin{pmatrix} 1&-1&-1\cr-1&1&-1\cr-1&-1&1
\end{pmatrix}$$
Since these Hessians are non-singular, we have a quadratic singularity
at each point, so each singular point contributes a single vanishing
cycle to the homology of other fibers of $\kappa$. We hence
have:
\begin{align*}
H_2(V_{2})&\cong \Z~\quad\text{generated by the vanishing cycle at
  }(0,0,0),\cr
H_2(V_{-2})&\cong \Z^4\quad\text{generated by the vanishing cycles at
  }(\pm2,\pm2,\pm2),\cr
H_2(V_{t})&\cong \Z^5\quad\text{otherwise, generated by all five vanishing cycles}.
\end{align*}
The vanishing cycle at a quadratic singularity is only defined up to
sign; a self-isomorphism of the singularity whose derivative has
negative determinant multiplies the cycle by $-1$.  We will choose a
basis of $\Z^5=\Z^4\oplus\Z$ by choosing vanishing cycles at
$(2,-2,-2)$ and $(0,0,0)$ and using the identity map between our
chosen local coordinates $\x,\y,\z$ at the four singular points
$(\pm1,\pm1,\pm1)$ to choose the vanishing cycles at the other three
singular points. We will order our basis according to the ordering of
singular points in (\ref{eq:cycles}) above.

\section{Action of $\Gamma$}\label{sec:action}

The Klein $4$-group $\Sigma$ of sign-changes permutes the first four
singular points (\ref{eq:cycles}) of $\kappa$ by the regular
permutation representation.  It acts as the identity map in the local
coordinates that we introduced above. It therefore simply permutes the
four vanishing cycles corresponding to these points.  Moreover, it
fixes the singular point $(0,0,0)$ and acts with positive determinant
in local coordinates there, so it fixes the corresponding vanishing
cycle. It hence acts on $\Z^5$ by the $4$-group whose non-trivial
elements are:
$$
\begin{pmatrix}
  0&1&0&0&0\cr
1&0&0&0&0\cr
0&0&0&1&0\cr
0&0&1&0&0\cr
0&0&0&0&1
\end{pmatrix}\,, \qquad
\begin{pmatrix}
  0&0&1&0&0\cr
0&0&0&1&0\cr
1&0&0&0&0\cr
0&1&0&0&0\cr
0&0&0&0&1
\end{pmatrix}\,,\qquad
\begin{pmatrix}
  0&0&0&1&0\cr
0&0&1&0&0\cr
0&1&0&0&0\cr
1&0&0&0&0\cr
0&0&0&0&1
\end{pmatrix}\,.
$$

The group $\PGL(2,\Z)$ is generated by $\PSL(2,\Z)$ and the involution
$\bigl({-1~0\atop \phantom{-}0~1}\bigr)$. We first discuss the group
$\PSL(2,\Z)$. It is the free product of cyclic groups of orders 2 and
3 generated by 
$$
\begin{pmatrix}0&-1\\1&0
\end{pmatrix}
\quad\text{and}\quad
\begin{pmatrix}1&-1\\1&0
\end{pmatrix}$$
respectively. These are
represented by automorphisms
\begin{align*}
\alpha\colon\quad&  X\mapsto Y,\quad Y\mapsto X^{-1}\\
\beta\colon\quad&  X\mapsto XY,\quad Y\mapsto X^{-1}\,.
\end{align*}

Recalling that an automorphism $\gamma$ of
$\pi$ acts on characters by $(\gamma\chi)(g)=\chi\circ\gamma^{-1}(g)$,
the action of the above two automorphisms on characters are:
\begin{align*}
\alpha_*\colon\quad(x,y,z)&\mapsto (y,x,xy-z)\\
\beta_*\colon\quad(x,y,z)&\mapsto (y,z,x)
\end{align*}
respectively.

The map $\alpha_*$ exchanges the first two singular points and fixes
the other three. Its derivative acts with determinant $1$ in the local
coordinates at the singular points. It thus acts on the homology
$\Z^5$ of a general fiber by:
$$
\begin{pmatrix}
0&1&0&0&0\cr
1&0&0&0&0\cr
0&0&1&0&0\cr
0&0&0&1&0\cr
0&0&0&0&1
\end{pmatrix}
$$

The map $\beta_*$ of $\C^3$ permutes the first three singular points
of (\ref{eq:cycles}) above and fixes the other two. It also acts with
determinant $1$ in the local coordinates at the singular points. It
thus acts on the homology $\Z^5$ of a general fiber by:
$$
\begin{pmatrix}
0&0&1&0&0\cr
1&0&0&0&0\cr
0&1&0&0&0\cr
0&0&0&1&0\cr
0&0&0&0&1
\end{pmatrix}
$$
The Klein $4$-group $\Sigma$ and the elements $\alpha$ and
$\beta$ generate the permutation representation of $S_4$ on the first
four coordinates of $\Z^5$, as promised in the theorem of the
Introduction.

Finally we discuss the involution 
\begin{equation*}
\begin{pmatrix} -1 & 0 \\ 0 & 1 \end{pmatrix} \in\PGL(2,\Z). 
\end{equation*}
The automorphism of $\pi$ defined by
$$\gamma\colon X\mapsto X^{-1},~Y\mapsto Y$$
represents this involution, and induces the map
$$\gamma_*\colon (x,y,z)\mapsto (x,y,xy-z)$$ 
on the character variety $\C^3$. 
This fixes all the singular points of $\kappa$ and its derivative acts with
determinant $-1$ in the local coordinates at each point. It thus
reverses the signs of all the vanishing cycles, so it acts on $\Z^5$
by the matrix
$$
\begin{pmatrix}
  -1&0&0&0&0\cr
0&-1&0&0&0\cr
0&0&-1&0&0\cr
0&0&0&-1&0\cr
0&0&0&0&-1
\end{pmatrix}\,.
$$

\section{Intersection form}\label{sec:intersection form}
Consider a general fiber $V_t$ (that is, $t\ne\pm2$). It is well known that
the vanishing cycle of a quadratic singularity has self-intersection
number $-2$. 
The vanishing cycles arising from the four
singularities on $V_2$ are clearly disjoint. Thus the intersection
form on $H_2(V_t)$ with respect to our basis is
$$
\begin{pmatrix}
    -2&0&0&0&a_1\\
0&-2&0&0&a_2\\
0&0&-2&0&a_3\\
0&0&0&-2&a_4\\
a_1&a_2&a_3&a_4&-2
\end{pmatrix},$$
where $a_1,\dots,a_4$ are still to be determined. Since 
an index 2 subgroup of $\Gamma$ permutes the first four basis
elements, $a_1=a_2=a_3=a_4$. Call this common value
$a$. We must show $a=\pm1$ (the sign depends on how we have oriented
the vanishing cycles over $t=2$ and $t=-2$ with respect to each other,
and is thus indeterminate).

The intersection matrix has determinant $32(1-a^2)$. 
Let $\hat V_t$ be the intersection of $V_t$ with a very large ball
around the origin and $M^3=\partial\hat V_t$ the ``link at
infinity'' of $V_t$. Then $V_t$ retracts to $\hat V_t$, so they
have the same homology and intersection form. The homology exact
sequence 
\begin{equation*}
H_2(\hat V_t) \to H_2(\hat V_t,M^3)\to H_1(M^3)\to 0 
\end{equation*}
and Poincar\'e-Lefschetz duality 
\begin{equation*}
H_2(\hat V_t,M^3)\cong H^2(\hat V_t) )\cong H_2(\hat V_t)^* 
\end{equation*}
yield a standard long exact sequence
$$
H_2(\hat V_t) \to H_2(\hat V_t)^*\to H_1(M^3)\to 0$$
where the
first arrow is intersection form. Thus  $H_1(M^3)$ is the
cokernel of the intersection matrix, so $a=\pm1$ if and only if
$H_1(M^3)$ is infinite (and otherwise $|H_1(M^3)|$ would equal
$32|1-a^2|)$). It thus remains to prove that $H_1(M^3)$ is indeed
infinite.

Consider $M^3$ as a regular neighborhood boundary of
the divisor at infinity of the closure $\overline V_t$ of $V_t$ in 
$\P^3(\C)=\C^3\cup\C\P^2$. Since the highest order term of $\kappa$ is
$xyz$, $\overline V_t$ intersects the projective plane at infinity in
the union of the three lines $x=0$, $y=0$, and $z=0$ (these are three
of the $27$ lines on a cubic surface; we will meet the other $24$
lines on $V_t$ in the next section). It is easy to check in local
coordinates that $\overline V_t$ is non-singular at infinity (and
hence globally non-singular except when $t=\pm 2$).  The three lines
at infinity of $V_t$ intersect in a cyclic configuration, so their
neighborhood boundary certainly has infinite first homology, as
desired.

In fact (see \cite{neumann1981}), the neighborhood boundary $M^3$ of
any cyclic configuration of rational curves in a non-singular complex
surface is a $T^2$ bundle over $S^1$ whose monodromy can be computed
from the self-intersection numbers $e_i$ of the curves as 
\begin{equation*}
A=\prod_i \begin{pmatrix} 0&-1\\1&-e_i \end{pmatrix}.  
\end{equation*}
In our case it is not hard to calculate that each of the three curves
has self-intersection $-1$ (this is true for all $27$ lines on a cubic
surface) so the monodromy is 
\begin{equation*}
\begin{pmatrix}
 0&-1\\1&1
\end{pmatrix}^3
=
\begin{pmatrix}
  -1&0\\0&-1
\end{pmatrix}
\end{equation*}
and hence 
\begin{equation*}
H_1(M^3)=\Z\oplus(\Z/2)^2. 
\end{equation*}
One can check that this
agrees with the cokernel of the intersection matrix of the
proposition.

\section{Topology of fibers via the 27 lines on a cubic}\label{sec:27 lines}

An alternative way of computing the homology of the fibers of
$\kappa$ and the action of $\Gamma$ involves the projection to 
a coordinate line, and provides further insight into the topology.
We sketch it briefly, hopefully with 
sufficient explanation that the reader can fill in the details. 
It is, in fact, the way we first obtained the results.

Suppose $t\ne\pm2$ and consider the projection $\pi\colon V_t \to \C$
given by the $z$-coordinate. The general fiber $\pi^{-1}(s)$ of this
projection is given by $x^2+y^2-xys=t+2-s^2$, which we can write as
$$(x-\lambda y)(x-\lambda' y)=t+2-s^2$$ 
with $\lambda,\lambda'=\frac12(s\pm\sqrt{s^2-4})$. Thus
$\pi^{-1}(s)\cong\C^*$ by the map $(x,y)\mapsto x-\lambda y$ unless
$\lambda=\lambda'$ or $t+2-s^2=0$. The general fiber of $\pi$ is thus
$\C^*$, while the special fibers occur for $s=\pm2$ and
$s=\pm\sqrt{t+2}$. These special fibers are:
\begin{align*}
\pi^{-1}(2)=L_1\cup L_2;\quad &L_1=\{x=y+\sqrt{t-2}\},\\&
L_2=\{x=y-\sqrt{t-2}\}\\
\pi^{-1}(\sqrt{t+2})=L_3\cup L_4;\quad &L_3=\hbox{$\{\frac
xy=\frac{\sqrt{t+2}+\sqrt{t-2}}2\}$},\\&
L_4=\hbox{$\{\frac xy=\frac{\sqrt{t+2}-\sqrt{t-2}}2\}$}\\
\pi^{-1}(-\sqrt{t+2})=L_5\cup L_6;\quad &L_5=\hbox{$\{\frac
xy=-\frac{\sqrt{t+2}+\sqrt{t-2}}2\}$},\\&
L_6=\hbox{$\{\frac xy=-\frac{\sqrt{t+2}-\sqrt{t-2}}2\}$}\\
\pi^{-1}(-2)=L_7\cup L_8;\quad &L_7=\{x=-y+\sqrt{t-2}\},\\&
L_8=\{x=-y-\sqrt{t-2}\}
\end{align*}
These $8$ lines $L_1,\dots,L_8$ and the
analogous lines coming from the projections to the $x$ and $y$
coordinates give $24$ lines, which with the $3$ lines at infinity
give the full $27$ lines on the projective cubic $\overline V_t$.

Choose small disks in $\C$ 
$$D_1, D_2, D_3, D_4\text{ around }z=2,
\sqrt{t+2}, -\sqrt{t+2}, -2,$$ 
and paths:
\begin{align*}
  \gamma_1&\text{ from $2$ to }\sqrt{t+2},\\
 \gamma_2&\text{ from $\sqrt{t+2}$ to }-\sqrt{t+2},\\
 \gamma_3&\text{ from $-\sqrt{t+2}$ to }-2,
\end{align*}
so that these paths do not intersect each other except at endpoints.

Then $X:=\pi^{-1}(D_1\cup\gamma_1\cup D_2\cup\gamma_2\cup
D_3\cup\gamma_3\cup D_4)$ is a deformation retract of $V_t$. Moreover
\begin{align*}
\pi^{-1}(D_1)\simeq S^2\simeq \pi^{-1}(D_4)&\qquad\text{(see below)}\\
\pi^{-1}(D_2)\simeq \hbox to 0pt{~~$*$\hss}\phantom{S^2}
\simeq \pi^{-1}(D_3)&\qquad\text{(since $\pi^{-1}(\pm\sqrt{t+2})\simeq
  *$)},\\
\end{align*}
so a simple Mayer-Vietoris calculation of homology of 
$X$ gives
$$\begin{matrix}
  H_2(V_t)\quad\cong&\Z&\oplus&\Z&\oplus&\Z&\oplus&\Z&\oplus&\Z\\
\text{generators:}&\alpha_1&&\alpha_2&&\alpha_3&&\alpha_4&&\alpha_5
\end{matrix}$$
where:
\begin{align*}
\alpha_1&\text{ is a generator of }H_2(\pi^{-1}(D_1),\\
\alpha_5&\text{ is a generator of }H_2(\pi^{-1}(D_4),\\
\alpha_2,\alpha_3,\alpha_4&\text{ are supported over
  }\gamma_1,\gamma_2,\gamma_3. 
\end{align*}
We shall see that the intersection form for this basis is
$$
\begin{pmatrix}
  -4&2&0&0&0\\
2&-2&-1&0&0\\
0&-1&-2&1&0\\
0&0&1&-2&2\\
0&0&0&2&-4
\end{pmatrix}$$

To describe the $\alpha_i$ more explicitly we draw the special fibers
and $\pi^{-1}(D_i)$ in the projective completion $\overline V_t$ of
$V_t$:
 $$
   \includegraphics[width=.8\hsize]{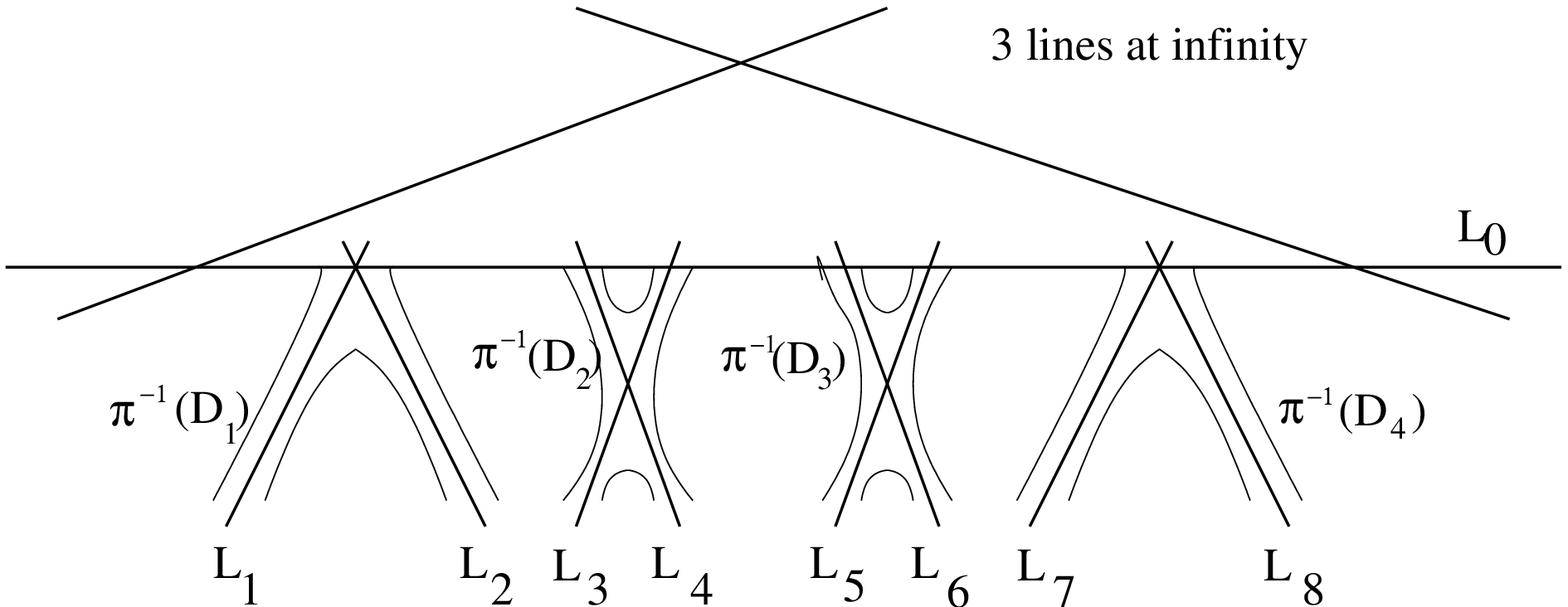}
 $$
Each line in the diagram has self-intersection number $-1$.
We first describe the classes $\alpha_2,\alpha_3,\alpha_4$.

For $\alpha_2$, choose an arc $\gamma$ in the line $L_0$ at infinity
from the point 
\begin{equation*}
L_0\cap L_2=[1:1] 
\end{equation*}
to the point 
\begin{equation*}
L_0\cap L_3=\bigg[{\frac{\sqrt{t+2}+\sqrt{t-2}}2}:1\bigg].  
\end{equation*}
We can form the
connected sum of $L_2$ and $-L_3$ by connecting them by a tube along
this arc (the tube is $s^{-1}(\gamma)$ where $s\colon\partial NL_0\to
L_0$ is the projection of the boundary of a tubular neighborhood of
$L_0$). This connected sum $L_2 \# -L_3$ lies in $V_t$ and 
represents the class $\alpha_2$.

We construct $\alpha_3$ and $\alpha_4$ similarly as $L_4 \# -L_5$ and
$L_5 \# -L_7$.  We shall see later that the choices of the tubes at
infinity to form these connected sums does not affect the homology
classes. 

The classes $\alpha_1$ and $\alpha_5$ are constructed similarly as
$L_1 \#-L_2$ and $L_7 \# -L_8$, but it is not so clear how the tubes are
chosen. To clarify this blow up the point $L_0\cap L_1$;
 $$
   \includegraphics[width=.6\hsize]{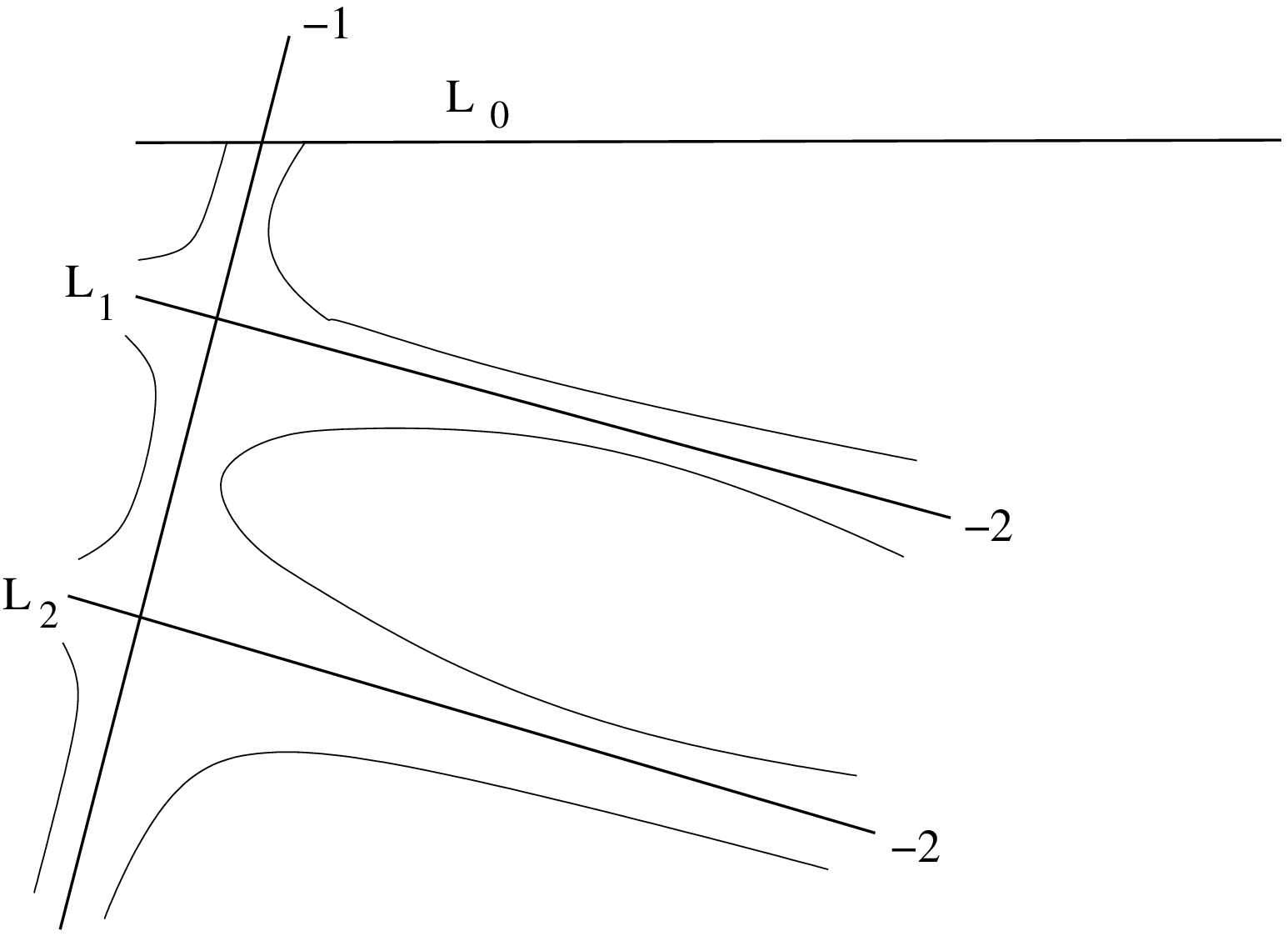}
 $$
then $\pi^{-1}(D_1)$ can be seen as the union of a portion
$N_E$ in a neighborhood of the exceptional curve $E$ and portions
$N_1$ and $N_2$ in neighborhoods of $L_1$, $L_2$.  $N_E$ is a
$(D^2-\{0\})$-bundle over $E-(E\cap L_0)$, hence homotopy equivalent
to a circle. Adding $N_1$ and $N_2$ to this just adds $2$-handles
along parallel copies of the circle, so $\pi^{-1}(D_1)\simeq S^2$, as
claimed earlier. The $S^2$ can be seen as $L_1 \#-L_2$ with the
connected sum formed along a tube over an arc in $E$. Its homology
class $\alpha_1$ thus has self-intersection number
\begin{equation*}
\alpha_1\cdot\alpha_1=L_1\cdot L_1+L_2\cdot L_2=-2+(-2)=-4 
\end{equation*}
(since blowing up on $L_1$ and $L_2$ has decreased their self-intersection
numbers by $1$). The next paragraph gives an alternative argument.

We can see mutual intersection numbers of all our elements as follows.
In the projective cubic $\overline V_t$ an element such as 
$\alpha_2 = L_2\# -L_3$ is homologous to the sum $L_2 + (-L_3)$, since the tube
used to perform connected sum bounds a solid tube. Since the inclusion
$V_t\to \overline V_t$ respects intersection number, we can simply
compute intersection numbers in $\overline V_t$. Thus, for instance,
$\alpha_2$ clearly has self-intersection number 
\begin{equation*}
L_2\cdot L_2 + L_3\cdot L_3 = -2. 
\end{equation*}

The claimed intersection matrix for the basis
$\alpha_1,\dots,\alpha_5$ is now easy to verify.

To see that the choices of the ``tubes at infinity'' do not affect the
homology classes $\alpha_2,\alpha_3,\alpha_4$, note that if we move
the path in $E_0$ defining such a tube past an intersection point
$E_0\cap E_i$ then we change the homology class by the class of a
torus $s^{-1}(\delta)$, where $\delta$ is a small loop in $E_0$ around
this point. Such a torus is supported entirely in $\pi^{-1}(D_i)$ and
has self-intersection number $0$, so its homology class is zero.

One can verify that the ordered basis of section \ref{sec:homology},
corresponding to the vanishing cycles at
\def\mb#1#2{\hbox to #1pt{\hss$#2$\hss}}
$$\mb{65}{(2,-2,-2),} \mb{65}{(-2,2,-2),}\mb{65}{(-2,-2,2),}
\mb{50}{(2,2,2),}\mb{45}{(0,0,0)},$$
is
$$\mb{65}{-(\alpha_4+\alpha_5),}\mb{65}{-\alpha_4,}
\mb{65}{\alpha_1+\alpha_2,}\mb{50}{\alpha_2,}\mb{45}
{-\alpha_3}.$$

\vspace{3pt}\noindent 
In terms of our lines these are:
$$\mb{65}{L_8 -L_5,}\mb{65}{L_7 -L_5,} \mb{60}{L_1-L_3,}
\mb{55}{ L_2-L_3,} \mb{50}{L_5-L_4.}$$

Recall that we only determined the basis of section \ref{sec:homology}
up to an overall sign. If we replace each $L_i$ by the other line in
the same fiber of $\pi$, so $L_8-L_5$ becomes $L_7-L_6$ for example,
each element of the basis is replaced by its negative.

The group $\Gamma$ acts on the projective cubic $\overline V_t$ by
rational maps. It has a linear subgroup isomorphic to $S_4$ (this
subgroup is generated by the group of sign changes and the elements
$\alpha_*\gamma_*$ and $\beta_*$ of section \ref{sec:action}). This
linear subgroup acts on $\overline V_t$ by morphisms (biholomorphic
maps) and therefore permutes the $27$ lines.  The whole of $\Gamma$ is
generated by this linear subgroup and the involution $\gamma_*$. This
latter involution acts on $\overline V_t$ by the rational map which
blows up the intersection point of the lines $x=0$ and $y=0$ at
infinity and blows down the line $z=0$ at infinity.  It is not hard to
derive the action on homology also from this point of view.

\section{The four-holed sphere}\label{sec:four holed sphere}

 Suppose $\psphere$ is a four-holed sphere. Its fundamental group $\pi$
 is freely generated by three peripheral elements $\delta_1,\delta_2,\delta_3$.
 The generator $\delta_4$ of the fourth boundary component satisfies
 the relation
 \begin{equation*}
 \delta_1 \delta_2 \delta_3 \delta_4 = 1
 \end{equation*}
 and we consider the redundant geometric presentation
 \begin{equation*}
 \langle X,Y,Z, \delta_1, \delta_2,\delta_3, \delta_4 \mid
 \delta_1 \delta_2 \delta_3 \delta_4 = 1, 
 X = \delta_1 \delta_2, 
 Y = \delta_2 \delta_3, 
 Z = \delta_3 \delta_1 \rangle.
 \end{equation*}
 Let $t = (t_1,t_2,t_3,t_4)\in\C^4$. 
 The relative character variety results by fixing the trace  
 \begin{equation*}
 \tr \rho(\delta_i) = t_i 
 \end{equation*}
 for $i=1,2,3,4$. In that case (see Goldman~\cite{Erg,fricke}, 
 Benedetto-Goldman~\cite{BenGol}, which is based on Magnus~\cite{Magnus}),
 the relative character variety is given by the cubic hypersurface
 in $\C^3$ defined by:
 \begin{align*}
 x^2 + y^2 + z^2 + x y z = 
 &(t_1t_2 + t_3t_4) x + 
 (t_1t_4 + t_2t_3) y + 
 (t_1t_3 + t_2t_4) z \\
 &+(4 - t_1^2 - t_2^2 - t_3^2 - t_4^2 - t_1 t_2 t_3 t_4).
 \end{align*}

Using the change of coordinates that reverses the signs of $x$, $y$
and $z$ puts this in the form:
 \begin{equation*}
 \kappa_{P,Q,R}(x,y,z)= S 
 \end{equation*}
where
$$\kappa_{P,Q,R}(x,y,z)=x^2 + y^2 + z^2 - x y z - P x - Q y - R z -
2$$
as in the Introduction, and $P=-(t_1t_2+t_3t_4)$, $Q=- (t_1t_4 +
t_2t_3) $, $R= - (t_1t_3 + t_2t_4)$, and $S= (2 - t_1^2 - t_2^2 - t_3^2
- t_4^2 - t_1 t_2 t_3 t_4)$.

\vspace{6pt} We first discuss the the group $\Aut(\kappa_{P,Q,R})$ of
polynomial automorphisms of the polynomial $\kappa_{P,Q,R}$.

Any affine automorphism $\phi=(\phi_1,\phi_2,\phi_3)$ of $\C^3$ that
preserves $\kappa_{P,Q,R}$ extends to $\P^3(\C)$ and must therefore
preserve the intersection of the closure $\overline V_t$ of the fiber
$V_t=\kappa_{P,Q,R}^{-1}(t)$ of $\kappa_{P,Q,R}$ with the plane at
infinity. This intersection consists of the three coordinate lines. It
follows that $\phi$ must permute the coordinates $x,y,z$ up to
multiplication by scalars. Inspecting $\kappa_{P,Q,R}$ we see that the
only non-trivial multiplication by scalars that can occur is
multiplying two coordinates by $-1$ if the corresponding two of
$P,Q,R$ are zero.  We thus see:
\begin{lemma}
  The group $L_{P,Q,R}=\Aut(\kappa_{P,Q,R})\cap\Aff(\C^3)$ is
  generated by permutation and sign-change automorphisms and is 
one of $$\{1\},\quad C_2,\quad C_2\times C_2,\quad S_3,\quad
S_4,$$ according as no two of $P,Q,R$ are equal; just two of them are
equal but nonzero; just two of them are equal to zero; all three
equal and nonzero; all three zero.\qed
\end{lemma}
Direct substitution shows that the following three involutions are in
$\Aut(\kappa_{P,Q,R})$:
\begin{align*}
  \tau_1\colon(x,y,z)&\mapsto (x,y,xy-z+R)\\
  \tau_2\colon  (x,y,z)&\mapsto (yz-x+P,y,z)\\
  \tau_3\colon  (x,y,z)&\mapsto (x,xz-y+Q,z)\,.
\end{align*}
\begin{theorem}
  $\Aut(\kappa_{P,Q,R})$ is generated by the above three involutions
  and $L_{P,Q,R}$. The three involutions generate a normal subgroup
  which is the free product
  $C_2*C_2*C_2=\langle\tau_1,\tau_2,\tau_3~|~
\tau_1^2=\tau_2^2=\tau_3^2=1\rangle$
  and
  $$\Aut(\kappa_{P,Q,R})= (C_2*C_2*C_2)\rtimes L_{P,Q,R}\,.$$
  (We
  identify this $C_2*C_2*C_2$ with $\PGL(2,\Z)_{(2)}$ in Theorem
  \ref{th:pgl}.)
\end{theorem}
\begin{proof}
  The proof follows Horowitz's proof \cite{horowitz} for the case $P=Q=R=0$.
  The group generated by the $\tau_i$ is normalized by $L_{P,Q,R}$.
  Let $\phi=(\phi_1,\phi_2,\phi_3)\in\Aut(\kappa_{P,Q,R})$ have degree
  $>1$. We claim there is a unique one of the three involutions
  $\tau_i$ for which $\tau_i\circ\phi$ has lower degree than $\phi$.
  It then follows that there is a unique reduced word $\omega$ in the
  $\tau_i$ such that $\omega\circ\phi\in L_{P,Q,R}$. This implies
  the theorem.
  
  To prove the claim, assume for convenience that
  $\deg\phi_1\le\deg\phi_2\le\deg\phi_3$ and $\deg\phi_3>1$. The
  highest order terms in the equation
  $$\phi_1^2+\phi_2^2+\phi_3^2-\phi_1\phi_2\phi_3-P\phi_1-Q\phi_2-R\phi_3=
  x^2 + y^2 + z^2 - x y z - P x - Q y - R z$$
  imply that $\deg
  \phi_3=\deg\phi_1+\deg\phi_2$ and that the highest order terms in
  $\phi_1\phi_2-\phi_3$ must cancel. It follows that $\tau_1\circ\phi$
  has lower degree than $\phi$. Clearly, neither of the other $\tau_i$
  has this property.
\end{proof}

\vspace{6pt} We now discuss the mapping class group
$\pi_0\homeo(\psphere,\partial_1 \psphere,\dots,\partial_4 \psphere)$
which fixes the boundary components of the four-holed sphere $N$.

We can array the punctures of $\psphere$ in order around the equator
of a $2$-sphere and place our basepoint on the part of the equator
joining the first two punctures.  Then the involution $T_1$ that
reflects north-south across the equator takes the elements
$X=\delta_1\delta_2$, $Y=\delta_2\delta_3$, and $Z=\delta_3\delta_1$
of $\pi_1\psphere$ to
$\delta_1^{-1}\delta_2^{-1}=\delta_1^{-1}X^{-1}\delta_1$,
$\delta_3^{-1}\delta_2^{-1}=Y^{-1}$, and $\delta_2\delta_4$. The
action on traces is $(x,y,z)\mapsto (x,y,xy-z+R)$. Similarly, by
arraying the four punctures in the orders $1324$ and $1243$ around the
equator we get two more involutions $T_2,T_3$; the actions on $\C^3$
of these three involutions are the three involutions
$\tau_1,\tau_2,\tau_3$ above.

It is not hard to see that the pairwise products $T_3T_1$, $T_1T_2$,
and $T_2T_3$ of these involutions are the Dehn twists $T_X,T_Y,T_Z$,
on the three separating simple closed curves that represent the
elements $X,Y,Z$ of $\pi_1\psphere$.  It is well-known that any two of
these Dehn twists freely generate $\pi_0\homeo^+(\psphere,\partial_1
\psphere,\dots,\partial_4 \psphere)$.  It follows that $T_1,T_2,T_3$
generate $\pi_0\homeo(\psphere,\partial_1 \psphere,\dots,\partial_4
\psphere)$.  Moreover, they generate it as a free product
$C_2*C_2*C_2$ since their action on $\C^3$ is as this free product.
Summarizing:

\begin{theorem}
  The mapping class group $\pi_0\homeo(\psphere,\partial_1
  \psphere,\dots,\partial_4 \psphere)$ is the free product
  $C_2*C_2*C_2$ generated by the three involutions $T_1,T_2,T_3$ above,
  and it acts on the relative character varieties via the action of 
$C_2*C_2*C_2$ of the previous theorem.\qed
\end{theorem}

Note that our computation implies that the Dehn twists $T_X$ and
$T_Y$ act
on the  relative character varieties by the following actions on
$\C^3$
\begin{align*}
  \tau_X\colon (x,y,z)&\mapsto (x,x^2y-xz+Rx-y+Q,xy-z+R)\\
  \tau_Y\colon  (x,y,z)&\mapsto (yz-x+P,y,y^2z-xy+Py-z+R)\,.
\end{align*}
This agrees with the calculation in \cite{Erg} (except for sign
differences because we changed coordinates to reverse the signs of
$x,y,z$).  But it should not agree, since our Dehn twists are actually the
inverses of the Dehn twists of \cite{Erg}.  The discrepancy is because
automorphisms $\gamma$ of $\pi_1\psphere$ act on characters $\chi$ by
$(\gamma,\chi)\mapsto\gamma\chi:=\chi\circ\gamma^{-1}$, but \cite{Erg}
uses $(\gamma,\chi)\mapsto\chi\circ\gamma$, which gives a
right-action.

\begin{theorem}\label{th:pgl}
  There are natural geometric isomorphisms:
  \begin{align*}
    \pi_0\homeo(\psphere,\partial_1
    \psphere,\partial\psphere-\partial_1\psphere)&\cong\PGL(2,\Z)\\
    \pi_0\homeo(\psphere,\partial_1
    \psphere,\dots,\partial_4 \psphere)&\cong\PGL(2,\Z)_{(2)}\\
    \pi_0\homeo(\psphere,\partial_\psphere)&\cong
  \PGL(2,\Z)_{(2)}\rtimes S_4\\
&\cong\PGL(2,\Z)\ltimes\Sigma\cong\Gamma
  \end{align*}
\end{theorem}
\begin{proof}
This is mostly well-known, so we describe it briefly.
$\GL(2,\Z)$ acts on $T^2=\R^2/\Z^2$ by linear maps, giving an
isomorphism of $\GL(2,\Z)$ with the mapping class group of $T^2$.
The central
element $-I\in\GL(2,\Z)$ acts as an involution fixing the four
half-integer points
$$\def\ss{\scriptstyle} (0,0) ,\quad({\ss\frac12},0),
\quad(0,{\ss\frac12}),\quad({\ss\frac12},{\ss\frac12}),$$
so the
quotient $\PGL(2,\Z)=\SL(2,\Z)/\{\pm I\}$ acts on the orbit space,
which is $S^2$, and permutes the four branch points.  This action
fixes the first branch point and permutes the other three via the
quotient $\PGL(\Z/2,2)$ (isomorphic to $S_3$).  Thus the kernel
$$\PGL(2,\Z)_{(2)}=\Ker(\PGL(\Z,2)\to\PGL(\Z/2,2))$$
acts on the
$2$-sphere fixing all four branch-points.  By replacing each branch
point by the circle of tangent directions at the point, we get actions
of $\PGL(2,\Z)$ and $\PGL(2,\Z)_{(2)}$ on the four-holed sphere
$\psphere$.  Since a mapping class of the four-holed sphere gives a
mapping class of the four-pointed sphere, which then lifts, up to
$C_2$ ambiguity, to a mapping class on the torus and thence to an
action on its homology, this gives the first two isomorphisms of the
theorem.

Finally, the mapping class
group $\pi_0\homeo(\psphere,\partial\psphere)$ acts on the set of boundary
components, giving a homomorphism to 
$S_4$, determining it as an $S_4$ extension of
$\PGL(2,\Z)_{(2)}=
\pi_0\homeo(\psphere,\partial_1\psphere,\dots,\partial_4\psphere)$:
$$
  1 \to \PGL(2,\Z)_{(2)}
\to 
\pi_0\homeo(\psphere,\partial\psphere) \to  S_4 \to 1
$$
This is a split extension, since we can map $S_4$ to
$\pi_0\homeo(\psphere,\partial\psphere)$ as follows: Construct $N$
from the boundary of a regular tetrahedron by removing neighborhoods
of the vertices; the tetrahedral group $S_4$ then acts by rigid
motions on $\psphere$. The Klein four-group $\Sigma\subset S_4\subset
\pi_0\homeo(\psphere,\partial\psphere)$ is the kernel of the action of
$\pi_0\homeo(\psphere,\partial\psphere)$ on homology of the torus
$T_2$ discussed above: each element of the Klein four-group comes from
a translation of $T^2=\R^2/\Z^2$ by an element of the half-integer
lattice, and these act trivially on homology.  Summarizing, we see
that $\pi_0\homeo(\psphere,\partial\psphere)$ fits in a diagram
$$\begin{CD}
  &&&&  1&&1\\
  &&&&@VVV @VVV\\
  &&&& \Sigma @>=>> \Sigma\\
  &&&&@VVV @VVV\\
  1 @>>> \PGL(2,\Z)_{(2)} @>>>\pi_0\homeo(\psphere,\partial\psphere)  
@>>> S_4 @>>> 1\\
  && @VV=V @VVV @VVV\\
  1 @>>> \PGL(2,\Z)_{(2)} @>>> \PGL(2,\Z)@>>> S_3 @>>> 1\\
  &&&&@VVV @VVV\\
  &&&&  1&&1\\
\end{CD}
$$
This is like the one for $\Gamma$ in the Introduction.  In fact,
both diagrams represent the group in the middle as the pullback for
the diagram
$$\begin{CD}
&& S_4 \\
&& @VVV\\
\PGL(2,\Z)@>>> S_3\\
\end{CD}
$$
so they are isomorphic, and
$\pi_0\homeo(\psphere,\partial\psphere)\cong\Gamma$.
\end{proof}
By explicitly
lifting the above generating involutions $T_1,T_2,T_3$ to the
torus and computing action on homology one finds that they correspond
under this isomorphism to the elements
$$
  \begin{pmatrix}
    1&0\\0&-1
  \end{pmatrix},\quad
  \begin{pmatrix}
    1&0\\2&-1
  \end{pmatrix},\quad
  \begin{pmatrix}
    1&2\\0&-1
  \end{pmatrix}\,
  $$
  respectively.  The three Dehn twists $T_X=T_3T_1$,
  $T_Y=T_1T_2$, $T_Z=T_2T_3$, therefore correspond
  to
$$\begin{pmatrix}
    1&2\\0&1
  \end{pmatrix},\quad
\begin{pmatrix}
    1&0\\-2&1
  \end{pmatrix},\quad
\begin{pmatrix}
    1&-2\\2&-3
  \end{pmatrix}\,.
$$

We now return to the topology of the fibers of the polynomials
$\kappa_{P,Q,R}$.  
\begin{theorem}
  Each fiber $V_t$ of $\kappa_{P,Q,R}$ has at most isolated
  singularities. It has reduced homology only in dimension $2$, and
  $H_2(V_t)=H_2(V)/I(V_t)$, where $V$ is a typical non-singular fiber
  and $I(V_t)$ is the ``subgroup of vanishing cycles'' for the
  singularities of $V_t$, which is non-trivial if and only if $V_t$ is
  singular. Moreover, $H_2(V)\cong\Z^5$, and it is the direct sum of
  the groups of vanishing cycles $I(V_t)$, summed over the singular
  fibers of $\kappa_{P,Q,R}$ (in particular, there are at most five
  singular fibers).
\end{theorem}
\begin{proof}
  Since the gradient $d\kappa_{P,Q,R}$ differs from the gradient
  $d\kappa$ by the constant vector $(P,Q,R)$, Lemma \ref{le:grad}
  applies to $\kappa_{P,Q,R}$, so $\kappa_{P,Q,R}$ is tame in
  Broughton's sense. Broughton \cite{broughton} (see also
  \cite{neumann-norbury}) proves that tameness implies the statements
  of the theorem, except, of course, for the fact that
  $H_2(V)\cong\Z^5$. We will see this below, where we show that each
  non-singular fiber $V$ of $\kappa_{P,Q,R}$ has the same topology as
  the non-singular fibers of $\kappa$.
\end{proof}
That $\kappa_{P,Q,R}$
has at most five singular fibers is also easy to see by direct
computation, and moreover, that for generic $P,Q,R$ it has exactly
five singular fibers, each with a single quadratic singularity.
However, more complex singularities can occur.  Benedetto and Goldman
discuss this topology in \cite{BenGol}, and describe, in particular,
under what conditions a relative character variety for $N$ can be a
singular fiber of the relevant cubic $\kappa_{P,Q,R}$.

In any case, we restrict now to a nonsingular fiber
$V=\kappa_{P,Q,R}^{-1}(t)$ of some $\kappa_{P,Q,R}$. Denote by
$\overline V$ the closure of $V$ in $\P^3(\C)$. As discussed in
sections \ref{sec:intersection form} and \ref{sec:27 lines} for the
case $P=Q=R=0$, $\overline V$ is a nonsingular projective cubic which
intersects the projective plane at infinity in the three lines $x=0$,
$y=0$, $z=0$, and these are three of the 27 lines on the cubic surface
$\overline V$, the others arising through singular fibers of the
projection of $V$ to the three coordinate axes. The computations are
essentially the same as the case $P=Q=R=0$. Thus the topology of the
nonsingular cubic $V$ does not change as the parameters $P,Q,R,t$ that
define it change. The action of the automorphism group of this cubic
on homology, which is always a subgroup of the group $\Gamma$, will
therefore also not change, so it is as described in the Introduction.


\begin{thebibliography}{99}

\bibitem{broughton} 
                Broughton, S. A., 
                Milnor number and the topology of polynomial hypersurfaces, 
                Inv. Math. {\bf92} (1988), 217--241.

\bibitem{clm1} 
                Cappell, S., Lee, R., Miller, E.Y.,
                The Torelli Group Action on Representation Spaces,
                Cont.\ Math.\ {\bf 258} (2000), 47--70

\bibitem{clm2} \bysame,
                The action of the Torelli group on the homology of 
                representation spaces is nontrivial, 
                Topology {\bf 39} (2000), 851--871.


\bibitem{BenGol} 
                Benedetto, R.\ and Goldman, W.,
                The topology of the relative character variety of the
                quadruply-punctured sphere,  
                Experimental Math. {\bf 8} 
                (1999), 85 - 103.

\bibitem{fricke} 
               Fricke, R, 
               \"Uber die Theorie der automorphen Modulgruppen, 
               Nachr. Akad. Wiss. G\"ottingen (1896), 91--101.


\bibitem{Erg} Goldman, W.M., 
                Ergodic Theory on Moduli Spaces, 
                Ann.\ Math.\ 146 (1997), 1-33

\bibitem{goldman} \bysame
               Action of the modular group on real $SL(2)$--characters of a 
               one-holed torus, Topology and Geometry (to appear).

\bibitem{goldman2} 
               \bysame, 
               An exposition of results of Fricke,
               (in preparation)

\bibitem{horowitz} 
               Horowitz, R.D. 
               Induced automorphisms on Fricke characters of free groups, 
               Trans. A.M.S. {\bf208} (1975), 41--50

\bibitem{Magnus} 
                Magnus, W.,
                {\em Rings of Fricke characters and automorphism groups
                of free groups,\/}
                Math.\ Zeit. {\bf 170} (1980), 91--103.

\bibitem{neumann1981} 
               Neumann, W.D., 
               A calculus for plumbing applied to the topology of complex 
               surface singularities and degenerating complex curves, 
               Trans. Amer.  Math. Soc. {\bf 268} (1981), 299--343.

\bibitem{neumann-norbury} \bysame and Norbury, P., 
               Unfolding singularities at infinity, 
               Math. Annalen {\bf318} (2000), 149--180.

\bibitem{nielsen} 
               Nielsen, J., 
               Die Isomorphismen der allgemeinen unendlichen Gruppe mit zwei 
               Erzeugenden, 
               Math. Ann. {\bf71} (1918), 385--397.

\end{thebibliography}
\end{document}